\documentclass[11pt]{article}
\usepackage{url}
\usepackage{amsmath}
\usepackage{amssymb}
\usepackage{amsfonts}

\title{Quadratic Frobenius probable prime tests costing two selfridges}
\author{Paul Underwood}
\date{\today}

\begin{document}
\maketitle

\begin{abstract}

By an elementary observation about the computation of the difference of squares for large integers,
 deterministic quadratic Frobenius probable prime tests are given with running times of approximately 2 selfridges.

\end{abstract}

\begin{section}{Introduction}

Much has been written about Fermat probable prime (PRP) tests~\cite{psw80,ra80,pk89}, Lucas PRP tests~\cite{ar97,wi98},
 Frobenius PRP tests~\cite{gr98,gr01,mu01,df03,se05,lo08,kh13}
  and combinations of these~\cite{bw80,po84,zh02}.
  These tests provide a probabilistic answer to the question: ``Is this integer prime?''
  Although an affirmative answer is not 100\% certain,
   it is answered fast and reliable enough for ``industrial'' use~\cite{cc99}.
    For speed, these various PRP tests are usually preceded by factoring methods such as sieving and trial division.

The speed of the PRP tests depends on how quickly multiplication and modular reduction can be computed during exponentiation.
 Techniques such as Karatsuba's algorithm~\cite[section 9.5.1]{cp05}, Toom-Cook multiplication,
  Fourier Transform algorithms~\cite[section 9.5.2]{cp05}
   and Montgomery exponentiation~\cite[section 9.2.1]{cp05}
    play their roles for different integer sizes. The sizes of the bases used are also critical.

Oliver Atkin introduced the concept of a ``Selfridge Unit''~\cite{at98}, approximately equal to the running
 time of a Fermat PRP test, which is called a {\em selfridge} in this paper.
   The Baillie-PSW test costs 1+3 selfridges,
   the use of which is very efficient when processing a candidate prime list.
    There is no known Baillie-PSW pseudoprime but Greene and Chen give a way to construct some similar counterexamples~\cite{gc99}.
     The software package Pari/GP implements a test similar to the Baillie-PSW test costing 1+2 selfridges.
      However, if the 2 selfridges Frobenius test presented in this paper is preceded with a Fermat 2-PRP test
      it also becomes 1+2 selfridges, but hitherto with the strength of a 1+1+2 selfridges test.

\end{section}

\begin{section}{Calculation$\pmod{n,x^2-ax+1}$}

At first sight, taking a modularly reduced power of $x$ would appear to be more efficient than
 taking a modularly reduced power of something more complicated, but this turns out to be false for the case presented here.

For most probable prime tests it is known that the given mathematical structures always work for primes
 and rarely work for composite numbers. Throughout this paper we are primarily concerned with the quotient ring $\mathbb{Z}_n[x]/(x^2-ax+1)$ where $f^{n+1} \equiv g \pmod{n, x^2-ax+1}$
 and $f$ and $g$ are polynomials in the polynomial ring $\mathbb{Z}_n[x]$, and $a$ is an integer.
 This means $f^{n+1}-g = nF + (x^2-ax+1)G$ where $F$ and $G$ are some polynomials.

 There is a recursive method to reduce any integer power of $x$ to a
 linear polynomial in $x$ since, for integer $k>1$, we have $x^k\equiv(ax-1)x^{k-2}$.
  For example, because $x^2 \equiv ax-1 \pmod{n,x^2-ax+1}$, then  $x^3  \equiv x(x^2)
         \equiv x(ax-1)
         \equiv ax^2-x
         \equiv a(ax-1) - x
         \equiv a^2x-a-x
         \equiv (a^2-1)x-a$.
However, we shall see below an iterative process which is by far superior to the recursive one.

For integer $a$, the equation $x^2-ax+1=0$ has discriminant $\Delta=a^2-4$ and solutions
$
x=\frac{a\pm\sqrt{\Delta}}{2}.
$
For an odd prime $p$, the Jacobi symbol $\Big(\frac{\Delta}{p}\Big)$ equals the Legendre symbol
 $\Big(\frac{\Delta}{p}\Big)$ and if $n$ is an odd integer in general then a negative Jacobi symbol
 implies a negative Legendre symbol, but the same is not guaranteed for positive Jacobi symbols.
  So if the Jacobi symbol $\Big(\frac{\Delta}{p}\Big)=-1$ then $\Delta$ will not be square modulo $p$,
   and by the Frobenius automorphism
\[
   x^p \equiv a-x \pmod{p,x^2-ax+1}
\]
so that
\[
\begin{array}{rcl}
x^p+x &\equiv &a\pmod{p,x^2-ax+1}\\
x^{p+1} &\equiv& 1 \pmod{p,x^2-ax+1}.
\end{array}
\]

In general, for a prime number, $p$, and for integers $S$ and $T$:
\[
(Sx+T)^p = S^px^p+(\sum^{p-1}_{i=1}\binom{p}{i}(Sx)^{p-i}T^i)+T^p
\]
and since the indicated binomial coefficients are divisible by $p$ and since
$S^p\equiv S\pmod{p}$ and $T^p\equiv T\!\!\pmod{p}$
we have
\[
(Sx+T)^p \equiv Sx^p+T \pmod{p}.
\]
Multiplying by $Sx+T$ gives
\[
\begin{array}{rcll}
(Sx+T)^{p+1} & \equiv & (Sx+T)(Sx^p+T)           & \pmod{p}\\
             & \equiv & S^2x^{p+1}+STx^p+STx+T^2 & \pmod{p}\\
             & \equiv & S^2x^{p+1}+ST(x^p+x)+T^2 & \pmod{p}\\
             & \equiv & S^2+aST+T^2              & \pmod{p,x^2-ax+1}.\;\;(*)
\end{array}
\]

In practice, left to right binary exponentiating of $Sx+T$ to the $(n+1)^{th}$ power can be accomplished with intermediate values $s$ and $t$ as follows. Firstly, obtain the binary representation of $n+1$. Secondly, assign $s=S$ and $t=T$. Thirdly, loop over bits of $n+1$, left to right, starting at the $2^{nd}$ bit, squaring the intermediate sum, $sx+t$, at each stage and if the corresponding bit is $1$ multiply the resulting squared intermediate sum by the base $Sx+T$.

Squaring the intermediate sum is achieved with appropriate modular reductions:
\[
\begin{array}{rcll}
(sx+t)^2 & =      & s^2x^2+2stx+t^2      &                  \\
         & \equiv & s^2(ax-1)+2stx+t^2   & \pmod{n,x^2-ax+1}\\
         & \equiv & (as^2+2st)x-s^2+t^2  & \pmod{n,x^2-ax+1}\\
         & \equiv & s(as+2t)x+(t-s)(t+s) & \pmod{n,x^2-ax+1}.
\end{array}
\]

If the bit is $1$ in the loop then the following must be calculated:
\[
\begin{array}{rcll}
(sx+t)(Sx+T)& =      & sSx^2+(sT+tS)x+tT    &                  \\
            & \equiv & sS(ax-1)+(sT+tS)x+tT & \pmod{n,x^2-ax+1}\\
            & \equiv & (asS+sT+tS)x+tT-sS   & \pmod{n,x^2-ax+1}.
\end{array}:
\]

If $a$, $S$ and $T$ are small then the squaring part is dominated by $2$ major multiplications and $2$ modular reductions:
$s$ by $as+2t$ modulo $n$, and $t-s$ by $t+s$ modulo $n$;
the ``if'' part is relatively faster. This makes an algorithm that is a little over $2$ selfridges.
The Pari/GP code for this process is
\begin{verbatim}
{general(n,a,S,T) = BIN=binary(n+1); LEN=length(BIN); aSpT=a*S+T; s=S; t=T;
 for(index=2, LEN, temp=(s*(a*s+2*t))%n; t=((t-s)*(t+s))%n; s=temp;
  if(BIN[index], temp=s*aSpT+t*S; t=t*T-s*S; s=temp));
   return( s==0 && t==(S*S+a*S*T+T*T)%n )}
\end{verbatim}

If $S=1$ and $T=0$ the program for computing the binary Lucas chain~\cite[algorithm 3.6.7]{cp05} is quicker,
 being $2$ selfridges:

\begin{verbatim}
{lucas_chain(n,a) = BIN=binary(n); LEN=length(BIN); va=2; vb=a;
 for(index=1, LEN,
  if(BIN[index], va=(va*vb-a)%n; vb=(vb*vb-2)%n, vb=(va*vb-a)%n; va=(va*va-2)%n));
   return( va==a && vb==2 )}
\end{verbatim}

\end{section}

\begin{section}{Equivalence of Tests}

The main test $(*)$ for an odd number $n$, presumed to be prime, with Jacobi symbol $\Big(\frac{\Delta}{n}\Big)=-1$, is
\[(Sx+T)^{n+1} \equiv S^2+aST+T^2 \pmod{n,x^2-ax+1}\]
and, by checking the discriminant, it is equivalent to
\[
y^{n+1} \equiv S^2+aST+T^2 \pmod{n,y^2-(aS+2T)y+S^2+aST+T^2}.
\]
Let $P = aS+2T$ and $Q = S^2+aST+T^2$ and the matrix
\[
A =
\left( \begin{array}{cc}
P&-Q\\
1&0
\end{array}\right).
\]

Note that, for our presumed prime $n$, $y^{n+1} \equiv Q \pmod{n, y^2-Py+Q}$ if and only if $A^{n+1} \equiv Q*I \pmod{n}$, where $I$ is the 2 by 2 identity matrix.
Using the multiplicative property of determinants for square matrices $X$ and $Y$ that $|X*Y| = |X||Y|$, we can ascertain that
  $|A^{n+1}| = |A|^{n+1} = Q^{n+1}$ and $|Q*I| = Q^2$. Thus $Q^{n-1} \equiv 1$ modulo $n$,
  since we assume $\gcd(PQ,n)=1$.

By Euler's criterion $Q^\frac{n-1}{2} \equiv \Big(\frac{Q}{n}\Big)$ modulo $n$ if $n$ is prime;
That is $n$ is Euler $Q$-PRP.
An implied 2 selfridges binary Lucas chain test can be shown to exist by letting $M = \frac{A^2}{Q*I}$ and then
\[
M^\frac{n+1}{2} \equiv \frac{A^{n+1}}{Q^\frac{n+1}{2}*I}\\
\equiv \frac{Q*I}{Q^\frac{n+1}{2}*I}\\
\equiv \frac{I}{Q^\frac{n-1}{2}*I}\\
\equiv \frac{I}{\Big(\frac{Q}{n}\Big)*I}\\
\equiv \Big(\frac{Q}{n}\Big)*I \pmod{n}.
\]
By the Cayley-Hamilton Theorem: Any 2 by 2 matrix $X$ satisfies its own characteristic equation $z^2 - trace(X)z + determinant(X) = 0$.
Given that
\[
M =
\left( \begin{array}{cc}
\frac{P^2}{Q}-1&-P\\
\frac{P}{Q}&-1
\end{array}\right)
\]
we can therefore deduce that
\[
z^\frac{n+1}{2}\equiv \Big(\frac{Q}{n}\Big) \pmod{n,z^2-(\frac{P^2}{Q}-2)z+1}.
\]

The number $21$ with $a=6$, $S=10$ and $T=4$, and so $P=68\equiv 5\pmod{21}$ and $Q=356\equiv 20\pmod{21}$,
is an example composite that passes the Euler PRP test but not the binary Lucas chain test.
 For a vice versa example: composite $27$ with $a=6$, $S=1$ and $T=7$, so that $P=20 \pmod{27}$ and $Q=92\equiv 11 \pmod{27}$,
  passes the binary Lucas chain test but not the Euler PRP test.

\end{section}

\begin{section}{The Main Algorithm for $S=1$ and $T=2$}

A test is now presented that is a little over $2$ selfridges. In comparison to the binary Lucas chain algorithm for a binary
 representation with an average number of ones and zeroes, the presented test requires an extra $7$ operations
  per loop iteration of multiple precision word additions or multiplications of multiple precision words by small numbers.
   Branching the code to handle the cases where $a=0$ and $a=1$ will reduce the extra operation count to $5$ and $6$ respectively.
    Perhaps the biggest difference in running times for the various PRP tests is that a Fermat PRP is dominated
     by a modularly reduced squaring per loop iteration;
      the binary Lucas chain test requires a modularly reduced squaring and a modularly reduced multiplication in its
       loop iteration; and the test presented in this section requires $2$ modularly reduced multiplications per loop iteration.
        As an example of this difference, if Fourier Transform arithmetic is used, only $1$
   forward transform is required for a squaring operation, whereas $2$ are required for multiplication.

For a candidate odd prime $n$, a {\em minimal} integer $a\geq0$ such that the Jacobi symbol $\Big(\frac{a^2-4}{n}\Big)=-1$ is sought. Then there is no ambiguity about how the algorithm works, there is no randomness. If, while searching for a minimum $a$, a value is found such that the Jacobi symbol $\Big(\frac{a^2-4}{n}\Big)=0$ then clearly $n$ is not prime, but this is unlikely to occur if sieving or trial division is performed firstly. Another reason for choosing a minimal $a$ is that there is more likelihood that the Jacobi symbol will be $0$ for the numerous candidates with small factors.
The time taken to test a Jacobi symbol is negligible, but some time can be saved by testing $a$ chosen
 in order from \[0,1,3,5,6,9,11,12,13,15,17,19,20,21,24,25,27,29,30,31,32\ldots\]
  Clearly, $2$ is to be omitted from this list. $a=4$ is omitted because it is covered by $a=0$ and $a=1$.
   $a=7$ is omitted since $\Big(\frac{3^2-4}{n}\Big)=\Big(\frac{7^2-4}{n}\Big)$, and so on.
Also, if a candidate prime equal to $1$ modulo $8$ is a square number then a Jacobi symbol equal to $-1$ will not be found.
 So it is recommended that a squareness test, which is rapid, is computed near the beginning.
  To ensure the implications of section 3,  the following is screened for:
\[
\gcd((a+4)(2a+5),n)=1.
\]
On finding a Jacobi symbol equal to $-1$ the following test can be done:
\[
(x+2)^{n+1} \equiv 2a+5 \pmod{n,x^2-ax+1}.
\]
The Pari/GP code for this process is
\begin{verbatim}
{selfridge2(n,a) = BIN=binary(n+1); LEN=length(BIN); ap2=a+2; s=1; t=2;
 for(index=2, LEN, temp=(s*(a*s+2*t))%n; t=((t-s)*(t+s))%n; s=temp;
  if(BIN[index], temp=ap2*s+t; t=2*t-s; s=temp));
   return( s==0 && t==(2*a+5)%n )}
\end{verbatim}

Using primesieve~\cite{ps14} and the GMP library~\cite{gmp14}, verification of the algorithm was pre-screened with the implied Fermat PRP test
$
(2a+5)^{n-1}\equiv 1 \pmod{n}.
$
For odd $n<2^{50}$ there were 1,518,678 such pseudoprimes.  The maximum $a$ required was $81$, for $n=170557004069761$. However, none that passed pre-screening were a pseudoprime for the full algorithm when run under Pari/GP~\cite{pari14}.

By examining the operations and their counts in the general test, it was decided by the author that the choice of $S=1$ and $T=2$ was optimal.
Moreover, one possible improvement in running times
might be achieved by using $S=1$ and $T=1$ for values of $2<a<n-2$, resulting in a hybrid test: Base $x+2$ could be used for
$a=0$ or $a=1$ and base $x+1$ used otherwise.
\end{section}
\pagebreak

\begin{section}{Conclusion}

We have seen how the algorithm in section 4 is effective for candidate odd primes less than $2^{50}$.
 No error rate bounds were examined but no failing pseudoprimes have been found.
Exploration of other $S$ and $T$ value pairs was not done.

 When implemented, a base $x+2$ Frobenius quadratic test costs about 2.5 selfridges which,
 when combined with a preceding one selfridge Fermat 2-PRP test, is not quite as fast as Pari/GP's ``ispseudoprime'' function, but it is faster by itself
  when testing a single sufficiently large number suspected of being prime.
  At the very large scale there is a handicap, especially when
   fast Fourier arithemetic comes into play because there are more forward transforms to be computed for using multiplications
    as opposed to using squaring operations.

    The Baillie-PSW test uses two independent tests:
     a strong Fermat 2-PRP test and a specific strong Lucas PRP test, whereas the test given in section 4 depends on one parameter,
      $a$.
    Can this difference influence reliability?

Figure 1 is a plot of pseudoprimes of the algorithm given in section 4 but for freely ranging $a$ and odd $n<2\cdot10^7$.
This leaves us with another question: ``Does a minimum $a$ for a pseudoprime ever come close to the minimum $a$
 required by the algorithm?'' A pseudoprime with a value of $a$ under say $n^\frac{1}{4}$ is extremely rare.
 There is no such $a$ for odd $n<2^{32}$.

\begin{figure}[!ht]
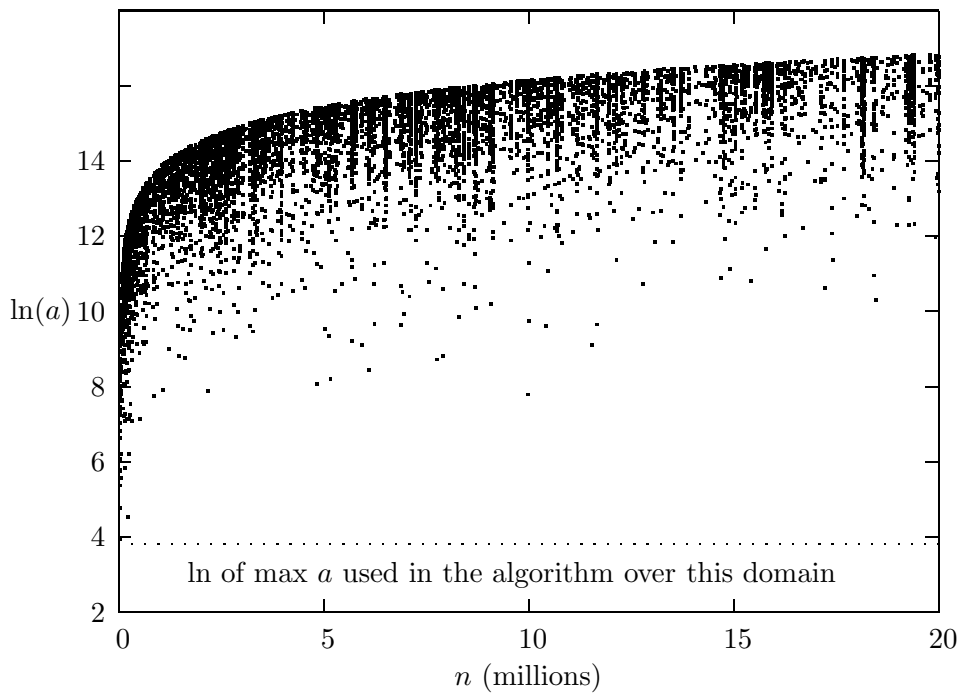

\centering
%\input{fu_plot.tex}
% GNUPLOT: LaTeX picture
\setlength{\unitlength}{0.240900pt}
\ifx\plotpoint\undefined\newsavebox{\plotpoint}\fi
\sbox{\plotpoint}{\rule[-0.200pt]{0.400pt}{0.400pt}}%

\caption{Pseudoprimes for S=1 and T=2}
\end{figure}

\end{section}

\begin{section}{Acknowledgements}

I thank Vincent Diepeveen for helping me code C/C++ and for allowing me to use his prime sieving function.
 Thanks too to members of mersenneforum.org and Yahoo! primenumbers groups for their encouraging remarks,
  in particular Carlos Pinho, Maximilian Hasler and Dana Jacobsen.

\end{section}

\bibliographystyle{ieeetr}
%\bibliography{fu}

\begin{thebibliography}{10}
\begin{footnotesize}
\bibitem{psw80}
C.~Pomerance, J.~L. Selfridge, and S.~S. Wagstaff, Jr., ``The pseudoprimes to
  $25 \cdot 10^9$'', {\em Mathematics of Computation}, vol.~35, no.~151,
  pp.~1003--1026, 1980.

\bibitem{ra80}
M.~O. Rabin, ``Probabilistic algorithm for testing primality'', {\em Journal of
  Number Theory}, vol.~12, no.~1, pp.~128--138, 1980.

\bibitem{pk89}
S.~H. Kim and C.~Pomerance, ``The probability that a random probable prime is
  composite'', {\em Mathematics of Computation}, vol.~53, no.~188,
  pp.~721--741, 1989.

\bibitem{ar97}
F.~Arnault, ``The {R}abin-{M}onier theorem for {L}ucas pseudoprimes'', {\em
  Mathematics of Computation}, vol.~66, no.~218, pp.~869--881, 1997.

\bibitem{wi98}
H.~C. Williams, {\em \'{E}douard {L}ucas and Primality Testing},
\newblock Wiley-Interscience, 1998.

\bibitem{gr98}
J.~Grantham, ``A {F}robenius probable prime test with high confidence'', {\em
  Journal of Number Theory}, vol.~72, pp.~32--47, 1998.

\bibitem{gr01}
J.~Grantham, ``Frobenius pseudoprimes'', {\em Mathematics of Computation},
  vol.~70, pp.~873--891, 2001.

\bibitem{mu01}
S.~M\"{u}ller, ``A probable prime test with very high confidence for n equiv 1
  mod 4'', {\em Proceedings of the 7th International Conference on the Theory
  and Application of Cryptology and Information Security: Advances in
  Cryptology}, pp.~87--106, 2001.

\bibitem{df03}
I.~B. Damg\r{a}rd and G.~S. Frandsen, ``An extended quadratic {F}robenius
  primality test with average and worst case error estimates'', {\em Lecture
  Notes in Computer Science. Fundamentals of Computation Theory (Springer
  Berlin Heidelberg)}, vol.~2751, pp.~118--131, 2003.

\bibitem{se05}
M.~Seysen, ``A {S}implified {Q}uadratic {F}robenius {P}rimality {T}est'',
  Cryptology ePrint Archive, Report 2005/462,
  \newblock \url{https://eprint.iacr.org/2005/462}

\bibitem{lo08}
D.~Loebenberger, ``A simple derivation for the {F}robenius pseudoprime test'',
  Cryptology ePrint Archive, Report 2008/124,
\newblock \url{https://eprint.iacr.org/2008/124}

\bibitem{kh13}
S.~Khashin, ``Counterexamples for {F}robenius primality test'', {\em eprint
  arXiv:1307.7920}, 2013.

\bibitem{bw80}
R.~Baillie and S.~S. Wagstaff, Jr., ``Lucas pseudoprimes'', {\em Mathematics of
  Computation}, vol.~35, pp.~1391--1417, October 1980.

\bibitem{po84}
C.~Pomerance, ``Are there counterexamples to the {B}aillie-{PSW} primality test?'',
  \url{http://www.pseudoprime.com/dopo.pdf}, 1984.

\bibitem{zh02}
Z.~Zhang, ``A one-parameter quadratic-base version of the {B}aillie-{PSW} probable
  prime test'', {\em Mathematics of Computation}, vol.~71, no.~240,
  pp.~1699--1734, 2002.

\bibitem{cc99}
C.~Caldwell, ``{P}robable {P}rime'',
  \url{http://primes.utm.edu/glossary/xpage/PRP.html}, 1999-2017.

\bibitem{cp05}
R.~Crandall and C.~Pomerance, {\em Prime Numbers, A Computational Perspective,
  2nd Ed.}
\newblock Springer, 2005.

\bibitem{at98}
A.~O.~L. Atkin., ``Intelligent primality test offer'', {\em Computational
  Perspectives on Number Theory (D. A. Buell and J. T. Teitelbaum, eds.),
  Proceedings of a Conference in Honor of A. O. L. Atkin, International Press},
  pp.~1--11, 1998.

\bibitem{gc99}
J.~R. Greene and Z.~Chen, ``Want to earn some cash the hard way?'',
  \url{http://www.d.umn.edu/~jgreene/baillie/Baillie-PSW.html}.

\bibitem{ps14}
``primesieve'', \url{http://primesieve.org}

\bibitem{gmp14}
``The {G}nu {M}ultiple {P}recision arithmetic library'', \url{https://gmplib.org}

\bibitem{pari14}
``Pari/{GP}'', \url{http://pari.math.u-bordeaux.fr}

\end{footnotesize}
\end{thebibliography}

\vspace{0.8cm}
\begin{footnotesize}
{\em E-mail address:} \texttt{paulunderwood@mindless.com}
\end{footnotesize}

\end{document}